\documentclass[10pt]{article}

\usepackage[french]{babel}
\usepackage{latexsym, amsmath, amssymb, graphicx, fancyhdr}

\newtheorem{lemme}{Lemme}[section]
\newtheorem{proposition}{Proposition}[section]

\newtheorem{theoreme}{Th\'eor\`eme}[section]
\newtheorem{corollaire}{Corollaire}[section]
\newtheorem{remarque}{\emph{Remarque}}[section]
\newtheorem{exemple}{\it{Exemple}}[section]

\title{{\bf Un mod\`{e}le non-lin\'eaire de prolif\'eration cellulaire : extinction des cellules et invariance}}
\author{Mostafa Adimy\thanks{M\'el : mostafa.adimy@univ-pau.fr} \quad et \quad Fabien Crauste\thanks{M\'el : fabien.crauste@univ-pau.fr}}
\date{Ann\'ee 2002}

\begin{document}
\thispagestyle{empty}

\maketitle

\begin{center}
{\em Laboratoire de Math\'ematiques Appliqu\'ees\\
Universit\'e de Pau et des Pays de l'Adour\\
Avenue de l'universit\'e, 64000 Pau, France}
\end{center}

\bigskip{}

\begin{abstract}
Dans ce travail, nous \'etudions un syst\`eme non-lin\'eaire structur\'e en \^age et maturit\'e d\'ecrivant la production des cellules sanguines dans la mo\"elle osseuse. Apr\`es int\'egration du mod\`ele, nous obtenons une \'equation aux d\'eriv\'ees partielles du premier ordre faisant appara\^\i tre un retard distribu\'e en temps et une d\'ependence non-locale dans la variable maturit\'e. Nous prouvons que l'unicit\'e des solutions ne d\'epend que des cellules de petites maturit\'es (cellules souches) et nous donnons un r\'esultat d'invariance.
\end{abstract}

\section{Introduction}

Nous nous int\'eressons dans ce travail \`a un mod\`ele
math\'ematique non-lin\'eaire d\'ecrivant une population
cellulaire : la production des cellules sanguines dans la mo\"elle
osseuse, \`a partir des cellules souches. Ce mod\`ele est d\'ecrit
par un syst\`eme de deux \'equations aux d\'eriv\'ees partielles
structur\'ees en \^age et maturit\'e.

Les mod\`eles de populations biologiques structur\'es en \^age et
maturit\'e sont apparus vers la fin des ann\'ees 60. Les premiers
\`a les avoir \'etudi\'es sont Keyfitz \cite{keyfitz} (1968),
Pollard \cite{pollard} (1973), Henry \cite{henry} (1976) et Mackey
\cite{mackey1978} (1978). Le mod\`ele de Mackey, dans lequel la
population cellulaire est suppos\'ee uniquement prolif\'erer, a
\'et\'e \'etudi\'e plus amplement \`a partir de 1990, notamment
num\'eriquement par Mackey et Rey dans \cite{mackey1992} (1992) et
\cite{mackey1993} (1993). En 1994, Mackey et Rudnicki
\cite{mackey1994} ont compl\'et\'e ce mod\`ele en tenant compte
d'une phase de repos en plus de la phase de prolif\'eration. Ils
obtiennent les \'equations sur lesquelles nous allons travailler
et qui sont pr\'esent\'ees en \ref{repos} et \ref{proliferation}.
Dans la phase de prolif\'eration, dont la dur\'ee est limit\'ee,
les cellules synth\'etisent de l'ADN puis se divisent, donnant
naissance \`a deux cellules filles. Chacune des cellules filles
entre ensuite imm\'ediatement dans la phase de repos. Elle peut y
demeurer toute sa vie et y mourir de vieillesse ou bien \^etre
introduite dans la phase de prolif\'eration \`a un certain moment.
Nous \'etudirons ici le cas non-lin\'eaire o\`u le taux de
r\'eintroduction d\'epend de la population totale de cellules au
repos.

Mackey et Rudnicki montrent que le syst\`eme d'\'equations aux
d\'eriv\'ees partielles qu'ils ont obtenu devient, apr\`es
int\'egration, un syst\`eme d'\'equations aux d\'eriv\'ees
partielles \`a retard. Dyson, Villella-Bressan et Webb en 1996
\cite{webb1996} et Adimy et Pujo-Menjouet en 2001 \cite{adimypujo}
ont montr\'e que la population de cellules de petites maturit\'es
influence le comportement de la population totale, dans le cas
o\`u les cellules sont suppos\'ees se diviser toujours au bout du
m\^eme temps.

Nous nous proposons d'\'etudier le mod\`ele propos\'e par Mackey et Rudnicki en 1994 \cite{mackey1994} en tenant compte d'une remarque faite par Mackey et Rey dans \cite{mackey1993} concernant la division cellulaire : nous supposerons que l'\^age de la cytocin\`ese (le point de division cellulaire) est distribu\'e selon une densit\'e de support $[\underline{\tau},\overline{\tau}]$, avec $0<\underline{\tau}<\overline{\tau}<+\infty$. Dans \cite{mackey1993}, Mackey et Rey n'ont \'etudi\'e que num\'eriquement ce mod\`ele et sans tenir compte de la phase de repos. Dans les travaux de Mackey et Rudnicki \cite{mackey1994}, Dyson, Villella-Bressan et Webb \cite{webb1996} et Adimy et Pujo-Menjouet \cite{adimypujo}, les cellules prolif\'erantes se divisaient toujours au bout du m\^eme temps fixe.

Enfin, rappelons que la production de cellules sanguines peut avoir deux comportements diff\'erents. Le premier correspond \`a une production normale de cellules et a lieu lorsque la population de cellules souches (c'est-\`a-dire les cellules de petites maturit\'es) est strictement positive. Le deuxi\`eme m\`ene \`a une production anormale de cellules et est le r\'esultat de la destruction ou de l'absence de cellules souches dans le milieu originel : ce cas correspond \`a l'an\'emie aplasique, une maladie s'attaquant aux cellules souches.

En supposant que la phase de prolif\'eration est suffisament longue, nous allons montrer que l'unicit\'e des solutions d\'epend uniquement des cellules souches et nous en d\'eduirons un r\'esultat d'invariance soulignant l'influence des cellules de petites maturit\'es sur le comportement de la population totale.

Notre travail est organis\'e de la fa\c con suivante : dans la section \ref{presentation}, nous pr\'esentons le mod\`ele biologique sur lequel est bas\'e notre travail ; dans la section \ref{section_equationsmodele}, nous obtenons, par int\'egration, le syst\`eme d'\'equations aux d\'eriv\'ees partielles d\'ecrivant l'\'evolution de la population ; dans la section \ref{para_formulationfaible}, nous donnons une formulation int\'egr\'ee du probl\`eme ; dans la section \ref{para_existenceunicite}, nous \'etablissons un r\'esultat d'unicit\'e des solutions ne d\'ependant que des cellules de petites maturit\'es et dans la section \ref{sectioninvariance} nous donnons un r\'esultat d'invariance.

\section{Pr\'esentation du mod\`ele} \label{presentation}

La vie d'une cellule d\'ebute avec la division d'une autre cellule, dite cellule m\`ere, et se termine, g\'en\'eralement, lors de sa propre division ; la cellule donne alors naissance \`a deux nouvelles cellules appel\'ees cellules filles. Durant sa vie, chaque cellule suit le m\^eme chemin que celui pris par sa m\`ere : il s'agit du cycle cellulaire. Le cycle cellulaire tel que nous l'entendons ici comprend deux phases : une phase de prolif\'eration et une phase de repos.

La phase de prolif\'eration est la phase active du d\'eveloppement cellulaire : elle est d\'edi\'ee \`a la synth\`ese de l'ADN et \`a la division cellulaire. Elle s'ach\`eve lors de la cytocyn\`ese : la division de la cellule proprement dite. Une cellule ne peut rester ind\'efiniment dans la phase de prolif\'eration : si elle ne se divise pas, elle meurt avec un taux $\gamma$.

Juste apr\`es leur naissance, les cellules filles entrent imm\'ediatement dans la phase de repos, ou phase $G_{0}$ (de l'anglais \emph{Gap}), qui est un stade quiescent du d\'eveloppement cellulaire : les cellules au repos sont inactives. Contrairement \`a ce qui se passe dans la phase de prolif\'eration, les cellules peuvent rester toute leur vie dans la phase de repos. La plupart seront introduites \`a un certain moment dans la phase de prolif\'eration et compl\`eteront ainsi le cycle cellulaire, mais les cellules au repos peuvent aussi mourir, avec un taux $\delta$.

Une cellule sera caract\'eris\'ee par son \^age et sa maturit\'e. La maturit\'e d\'esigne un aspect morphologique du d\'eveloppe\-ment de la cellule, au m\^eme titre que la taille ou l'\^age ; il s'agit d'une variable continue : dans l'\'evaluation de la maturit\'e, on tient compte, par exemple, du niveau de synth\`ese d'ADN, d'ARN, de la pr\'esence de mitochondries, etc. La d\'efinition de la maturit\'e est li\'ee \`a la nature de la cellule. Nous supposerons que la maturit\'e $m$ d'une cellule varie de $m=0$ \`a $m=1$.

Toutes les cellules proviennent des cellules embryonnaires issues des premiers jours du d\'eveloppement, appel\'ees cellules souches (\emph{stem cells}). Il s'agit de cellules indiff\'erenci\'ees et immatures (c'est-\`a-dire de petites maturit\'es) et capables, d'une part, de se multiplier (et donc d'augmenter la quantit\'e de cellules souches), et d'autre part, de donner naissance \`a des cellules diff\'erenci\'ees (cellules du foie, cellules musculaires, cellules sanguines, \ldots).

Il existe quatre types de cellules souches dont les cellules souches h\'emato\-po\"\i \'etiques de la mo\"elle osseuse. Ce sont celles-ci qui nous int\'eressent car elles sont \`a l'origine de toutes les cellules sanguines : globules rouges, blancs, plaquettes, etc.

Mis \`a part leur importance dans le d\'eveloppement d'une population cellulaire, il faut noter que les cellules souches sont aussi \`a l'origine de certaines maladies, graves, dont l'an\'emie aplasique. Il s'agit d'une maladie rare et s\'erieuse caract\'eris\'ee par l'absence de formation ou de division des cellules souches h\'emato\-po\"\i \'etiques, entra\^\i nant une chute de la production des cellules sanguines.

Plusieurs facteurs expliquent l'apparition de cette maladie.  Dans
50 \`a 65 \% des cas la nature de la maladie est idiopathique ou
primitive, c'est-\`a-dire qu'elle ne comporte pas de cause connue.
Dans les autres cas, elle est h\'er\'editaire ou due notamment \`a
des expositions aux radiations nucl\'eaires, \`a des toxines de
l'environnement (insecticides, vapeurs de benz\`ene, \ldots) ou
bien au traitement du cancer par radioth\'erapie, chimioth\'erapie
ou d'autres m\'edicaments.

\subsection{La phase de repos} \label{repos}

Nous d\'esignons par $n(t,m,a)$ la densit\'e de cellules dans la phase de repos \`a l'instant $t$, de maturit\'e $m$ et d'\^age $a$. Nous notons
\begin{displaymath}
    N(t,m)=\int_{0}^{\infty} n(t,m,a) da,
\end{displaymath}
le nombre total de cellules de maturit\'e $m$ dans la phase de repos \`a l'instant $t$.

La densit\'e de cellules $n(t,m,a)$ v\'erifie l'\'equation de conservation :
\begin{equation} \label{equation_n(t,m,a)}
\frac{\partial n}{\partial t} + \frac{\partial n}{\partial a} + \frac{\partial (Vn)}{\partial m} = -(\delta +\beta )n,
\end{equation}
o\`u $\delta$ est le taux de mortalit\'e, $\beta$ le taux de r\'eintroduction des cellules dans la phase de prolif\'eration et $V : [0,1] \rightarrow [0, +\infty )$ est la vitesse de maturation des cellules.

$\delta=\delta(m)$ et $\beta=\beta(m,N(t,m))$ sont suppos\'ees positives et continues, et $\beta$ est suppos\'ee d\'ecroissante par rapport \`a $N(t,m)$. $V$ est de classe $C^{1}$ sur $[0,1]$ et v\'erifie $V(0)=0$ et $V(m) > 0$ pour tout $m \in (0,1]$. De plus, nous supposerons que
\begin{displaymath}
\int_{0}^{m} \frac{ds}{V(s)} = + \infty,  \quad \textrm{pour tout} \  m \in (0,1].
\end{displaymath}
Cela signifie qu'une cellule de maturit\'e nulle ne peut jamais atteindre un \'etat o\`u sa maturit\'e serait non nulle.

\begin{exemple}
Si $V(s)\underset{s\to 0}{\sim}\alpha s^{p}$, pour tout $s\in[0,1]$, avec $\alpha>0$ et $p\geq 1$,  alors $\int_{0}^{m} \frac{ds}{V(s)}$ est une int\'egrale de Riemann divergente pour tout $m \in (0,1]$.
\end{exemple}

$n$ satisfait de plus la condition initiale :
\begin{equation} \label{condition_initiale_n(t,m,a)}
n(0,m,a)=\Upsilon (m,a), \quad \textrm{pour tout} \ (m,a)\in [0,1]\times [0,+\infty),
\end{equation}
o\`u $\Upsilon$ est une fonction suffisament r\'eguli\`ere v\'erifiant $\lim_{a \to \infty} \Upsilon (m,a) = 0$.

        \subsection{La phase de prolif\'eration} \label{proliferation}

Dans la phase de prolif\'eration, l'\^age d'une cellule varie de $0$ \`a $\overline{\tau}$ o\`u  $\overline{\tau}$ est la dur\'ee maximale de cette phase : parvenue \`a l'\^age  $\overline{\tau}$ une cellule doit s'\^etre divis\'ee sinon elle meurt. Cela entra\^ine qu'une cellule ne peut rester ind\'efiniment dans la phase de prolif\'eration.

Soit $p(t,m,a)$ la densit\'e de cellules de maturit\'e $m$ et d'\^age $a$ dans la phase de prolif\'eration \`a l'instant $t$. Nous notons
\begin{displaymath}
P(t,m)=\int_{0}^{\overline{\tau}} p(t,m,a) da,
\end{displaymath}
le nombre total de cellules de maturit\'e $m$ dans la phase de prolif\'eration \`a l'instant $t$.

Soit $\gamma$ une fonction de la maturit\'e $m$, positive et continue, repr\'esentant le taux de mortalit\'e des cellules dans la phase de prolif\'eration.

$p(t,m,a)$ v\'erifie l'\'equation de conservation :
\begin{equation} \label{equation_p(t,m,a)}
\frac{\partial p}{\partial t} + \frac{\partial p}{\partial a} + \frac{\partial (Vp)}{\partial m} = -\gamma p,
\end{equation}
o\`u $V$ repr\'esente la vitesse de maturation des cellules dans cette phase, d\'efinie en \ref{repos}. Nous supposons que la vitesse de maturation des cellules est la m\^eme dans les deux phases.

$p$ v\'erifie en plus la condition initiale :
\begin{equation} \label{condition_initiale_p(t,m,a)}
p(0,m,a)=\Gamma (m,a), \quad \textrm{pour tout} \ (m,a)\in [0,1]\times [0,\overline{\tau}],
\end{equation}
o\`u $\Gamma$ est une fonction suffisament r\'eguli\`ere, d\'efinie sur $[0,1]\times [0,\overline{\tau}]$.
\bigskip{}

Pour \^etre bien pos\'e, le syst\`eme d'\'equations (\ref{equation_n(t,m,a)})-(\ref{equation_p(t,m,a)}) doit \^etre muni de conditions aux bords, traduisant le comportement des cellules lors des changements de phases.

\subsection{Conditions aux bords}

En accord avec les remarques de Mackey et Rey (\cite{mackey1993}), nous supposerons que le point de cytocin\`ese des cellules prolif\'erantes est distribu\'e sur un intervalle $[\underline{\tau},\overline{\tau}]$, o\`u $0<\underline{\tau}<\overline{\tau}<+\infty$, selon une densit\'e $a\mapsto k(m,a)$. $\overline{\tau}$ repr\'esente la dur\'ee maximale de la phase de prolif\'eration tandis que $\underline{\tau}$ est choisi strictement positif car les cellules ne se divisent imm\'ediatement apr\`es \^etre entr\'ees en phase de prolif\'eration que dans le cas de certaines maladies qui ne nous int\'eressent pas ici.

Nous notons $g(m)$ la maturit\'e d'une cellule fille \`a sa naissance, lorsque la maturit\'e de la cellule m\`ere valait $m$. La fonction $g:[0,1]\to[0,1]$ est suppos\'ee continue, strictement croissante sur $[0,1]$, de classe $\mathcal{C}^{1}$ sur $[0,1)$ et telle que $g(m)<m$ pour $m \in (0,1)$. $g$ est alors inversible et nous prolongeons $g^{-1}$ en posant $g^{-1}(m)=1$ si $m>g(1)$.

L'application $k$ est telle que $m \mapsto k(m,a)$ est la proportion de cellules m\`eres d'\^age $a$ dans la phase de prolif\'eration pouvant donner apr\`es division une cellule fille de maturit\'e $m$.  $k$ est suppos\'ee positive et continue sur $[0,1]\times [\underline{\tau},\overline{\tau}]$, et telle que $k(m,a)=0$ si $m\geq g(1)$. Ceci par d\'efinition de la fonction $g$ : la maturit\'e d'une cellule m\`ere ne pouvant \^etre sup\'erieure \`a $1$, celle d'une cellule fille ne peut exc\'eder la valeur $g(1)$. Nous consid\`ererons donc, dans toute la suite, que la maturit\'e d'une cellule fille est comprise entre $0$ et $g(1)$.

La premi\`ere condition au bord est alors donn\'ee par
\begin{equation} \label{condbordn}
n(t,m,0)=2 \int_{\underline{\tau}}^{\overline{\tau}} k(m,a)p(t,g^{-1}(m),a)da.
\end{equation}

Elle signifie qu'une cellule de maturit\'e $m$ peut provenir d'un certain nombre de cellules m\`eres, d'\^ages diff\'erents, contrairement aux mod\`eles pr\'ecedemment \'etudi\'es par Mackey et Rudnicki \cite{mackey1994}, Dyson, Villella-Bressan et Webb \cite{webb1996} et Adimy et Pujo-Menjouet \cite{adimypujo}, o\`u la division survenait toujours au bout d'un temps fixe (la fin de la phase de prolif\'eration).

La seconde condition aux limites est
\begin{equation} \label{condbordp_nonlineaire}
p(t,m,0)=\int_{0}^{\infty} \beta (m,N(t,m))n(t,m,a)da = \beta (m,N(t,m)) N(t,m).
\end{equation}
Elle repr\'esente le flux de la population de cellules, de la
phase de repos vers la phase de prolif\'eration.

Nous consid\'erons le flot $\pi_{s} : [0,1] \to [0,1]$, d\'efini pour $s \leq 0$, solution de l'\'equation diff\'erentielle ordinaire:
\begin{displaymath}
\left\{ \begin{array}{rcll}
\displaystyle\frac{du}{ds}(s) & = & V(u(s)),& s \leq 0, \\
u(0)         & = & m, &
\end{array} \right.
\end{displaymath}
qui repr\'esente l'\'evolution de la maturit\'e des cellules pour atteindre une maturit\'e $m$ \`a l'instant $0$ \`a partir d'un temps $s \leq 0$. $\pi_{s}$ v\'erifie $\pi_{0}(m)=m$, $\pi_{s}(0)=0$ et $\pi_{s}(m) \in (0,1]$ pour $s\leq 0$ et $m \in (0,1]$.

En fait, nous pouvons remarquer (cf. \cite{pujo}) que $\pi_{s}$ peut s'\'ecrire explicitement
\begin{displaymath}
\pi_{s}(m)=h^{-1}(h(m)e^{s}), \quad \textrm{pour } m \in [0,1] \textrm{ et } s \leq 0,
\end{displaymath}
o\`u $h :[0,1] \to [0,1]$ est d\'efini par
\begin{displaymath}
h(m)=\left\{ \begin{array}{ll}
exp\bigg(-\displaystyle\int_{m}^{1} \frac{ds}{V(s)} \bigg),& \quad \textrm{pour } m \in (0,1],\\
0, & \quad \textrm{pour } m=0.
\end{array} \right.
\end{displaymath}
Notons que $h$ est une fonction continue et strictement croissante, donc inversible.

\begin{exemple}
Si $V(m)=\alpha m$, pour tout $m\in[0,1]$, avec $\alpha>0$, alors
\begin{displaymath}
h(m)=m^{1/\alpha}, \quad \forall m\in[0,1].
\end{displaymath}
Donc, pour tout $m\in[0,1]$,
\begin{displaymath}
h^{-1}(m)=m^{\alpha},
\end{displaymath}
et
\begin{displaymath}
\pi_{s}(m)=me^{\alpha s}, \quad \textrm{ pour } s\leq 0.
\end{displaymath}
Nous obtenons les fonctions utilis\'ees par Rey et Mackey dans \cite{mackey1992} et \cite{mackey1993} et Dyson, Villella-Bressan et Webb dans \cite{webb1996}.
\end{exemple}


\section{Equations du mod\`ele} \label{section_equationsmodele}

Avant de d\'eterminer les \'equations du mod\`ele sur lesquelles nous allons travailler nous \'etablissons un lemme qui nous sera utile.

Soit $\Delta$ l'application d\'efinie de $[0,+\infty)\times[0,g(1)]$ dans $\mathbb{R}^{+}$ par
\begin{displaymath}
\Delta(s,m) = \pi_{-s}(g^{-1}(m))=h^{-1}(h(g^{-1}(m))e^{-s}).
\end{displaymath}

\begin{lemme} \label{fonctiondelta}
$\Delta$ satisfaits les assertions suivantes :
\begin{itemize}
\item[(i)] l'application $s \in [0,+\infty) \mapsto \Delta(s,m)$ est d\'ecroissante pour chaque $m \in [0,g(1)]$,
\item[(ii)] l'application $m \in[0,g(1)]\mapsto \Delta(s,m)$ est croissante pour chaque $s\geq 0$, donc, pour tout $s\geq 0$ et tout $m \in [0,g(1)]$,
\begin{displaymath}
0=\Delta(s,0) \leq \Delta(s,m) \leq \Delta(s,g(1))=h^{-1}(e^{-s}),
\end{displaymath}
\item[(iii)] pour tout $m \in (0,g(1)]$,
\begin{displaymath}
\Delta(s,m) < m  \quad \textrm{ssi} \quad s > \int_{m}^{g^{-1}(m)} \frac{d\theta}{V(\theta)},
\end{displaymath}
\item[(iv)] pour tous $\sigma, s \geq 0$, et tout $m \in[0,g(1)]$,
\begin{displaymath}
\pi_{-\sigma}(\Delta(s,m))=\Delta(s+\sigma,m).
\end{displaymath}
\end{itemize}
\end{lemme}
\vspace{2ex}
\textbf{Preuve :} Les points \emph{(i)} et \emph{(ii)} sont \'evidents par d\'efinition, les fonctions $h$ et $g$ \'etant croissantes sur $[0,1]$.

Pour montrer le point \emph{(iii)}, remarquons que $\displaystyle\int_{m}^{g^{-1}(m)} \frac{d\theta}{V(\theta)} = \ln \frac{h(g^{-1}(m))}{h(m)}$ pour tout $m>0$. Alors nous avons, pour tout $m \in (0,g(1)]$ et tout $s\geq 0$,
\begin{eqnarray}
\Delta(s,m) < m & \Leftrightarrow & h(g^{-1}(m))e^{-s} < h(m) \nonumber \\
        & \Leftrightarrow & -s < \ln \frac{h(m)}{h(g^{-1}(m))} \nonumber \\
        & \Leftrightarrow & s > \int_{m}^{g^{-1}(m)} \frac{d\theta}{V(\theta)}. \nonumber
\end{eqnarray}
D'o\`u le r\'esultat.

Montrons \emph{(iv)}. Par d\'efinition nous avons, pour tous $\sigma, s \geq 0$ et tout $m \in[0,g(1)]$,
\begin{eqnarray}
\pi_{-\sigma}(\Delta(s,m)) & = & h^{-1}(h(h^{-1}(h(g^{-1}(m))e^{-s}))e^{-\sigma}) \nonumber \\
               & = & h^{-1}(h(g^{-1}(m))e^{-(s+\sigma)}) \nonumber \\
               & = & \Delta(s+\sigma,m). \nonumber
\end{eqnarray}
Le lemme est ainsi d\'emontr\'e. $\Box$\bigskip{}

Nous utiliserons, dans toute la suite, les notations suivantes :
\begin{displaymath}
    \xi (m,t)=exp\bigg\{-\int_{0}^{t} \gamma (\pi _{-s}(m))+V'(\pi _{-s}(m))ds\bigg\},
\end{displaymath}
pour tout $ m \in [0, 1]$ et tout $t \geq 0$, et
\begin{displaymath}
\zeta(m,a)=k(m,a)\xi(g^{-1}(m),a),
\end{displaymath}
pour tout $ m \in [0, 1]$ et tout $a \geq 0$.

Nous \'etablissons le r\'esultat suivant.

\begin{proposition} \label{equationsmodele}
Le nombre total, $N(t,m)$, de cellules de maturit\'e $m$ dans la phase de repos, v\'erifie les \'equations suivantes :

\vspace{1ex}
(i) Si $0 \leq t \leq \underline{\tau}$ et $ m \in [0, g(1)]$,
\begin{displaymath}
\frac{\partial }{\partial t}N(t,m)+\frac{\partial }{\partial m} (V(m)N(t,m)) = -(\delta (m)+\beta (m,N(t,m)))N(t,m)
\end{displaymath}
\begin{equation} \label{N(t,m)_1bis}
+2\xi (g^{-1} (m),t)\int_{\underline{\tau}}^{\overline{\tau}} k(m,a)\Gamma (\Delta(t,m), a-t) da.
\end{equation}

\vspace{1ex}
(ii) Si $\underline{\tau} \leq t \leq \overline{\tau}$ et $ m \in [0, g(1)]$,
\begin{displaymath}
\frac{\partial }{\partial t}N(t,m)+\frac{\partial }{\partial m} (V(m)N(t,m)) = -(\delta (m)+\beta (m,N(t,m)))N(t,m)
\end{displaymath}
\begin{displaymath}
\ +2\int_{\underline{\tau}}^{t} \zeta (m,a) \beta (\Delta(a,m),N(t-a,\Delta(a,m)))N(t-a,\Delta(a,m))da
\end{displaymath}
\begin{equation} \label{N(t,m)_1}
+2\xi (g^{-1} (m),t)\int_{t}^{\overline{\tau}} k(m,a)\Gamma (\Delta(t,m), a-t) da.
\end{equation}

\vspace{1ex}
(iii) Si $t \geq \overline{\tau}$ et $ m \in [0, g(1)]$,
\begin{displaymath}
\frac{\partial }{\partial t}N(t,m)+\frac{\partial }{\partial m} (V(m)N(t,m)) = -(\delta (m)+\beta (m,N(t,m)))N(t,m)
\end{displaymath}
\begin{equation}\label{N(t,m)_2}
+2\int_{\underline{\tau}}^{\overline{\tau}} \zeta (m,a)\beta (\Delta(a,m),N(t-a,\Delta(a,m)))N(t-a,\Delta(a,m))da.
\end{equation}

\vspace{1ex}
(iv) Si $t\geq 0$ et $m \in [g(1),1]$,
\begin{equation} \label{N(t,m)_3}
\frac{\partial }{\partial t}N(t,m)+\frac{\partial }{\partial m} (V(m)N(t,m)) = -(\delta (m)+\beta (m,N(t,m)))N(t,m).
\end{equation}

De m\^eme, la population totale, $P(t,m)$, de cellules de maturit\'e $m$ dans la phase de prolif\'eration v\'erifie les \'equations suivantes :

\vspace{1ex}
(v)  Si $0 \leq t \leq \overline{\tau}$ et $m \in [0,1]$,
\begin{eqnarray} \label{P(t,m)_1}
\frac{\partial }{\partial t} P(t,m)+\frac{\partial }{\partial m} (V(m)P(t,m)) & = & -\gamma (m)P(t,m)+\beta (m,N(t,m))N(t,m) \nonumber \\
    & - & \Gamma (\pi_{-t}(m),\overline{\tau} -t) \xi (m,t).
\end{eqnarray}

\vspace{1ex}
(vi) Si $t \geq \overline{\tau}$ et $m \in [0,1]$,
\begin{displaymath}
\frac{\partial }{\partial t} P(t,m)+\frac{\partial }{\partial m} (V(m)P(t,m)) = -\gamma (m)P(t,m)+\beta (m,N(t,m))N(t,m)
\end{displaymath}
\begin{equation}\label{P(t,m)_2}
\qquad \qquad\qquad -\xi (m,\overline{\tau})\beta (\pi_{-\overline{\tau}}(m),N(t-\overline{\tau},\pi_{-\overline{\tau}}(m)))N(t-\overline{\tau},\pi_{-\overline{\tau}}(m)).
\end{equation}
\end{proposition}
\vspace{2ex}
\textbf{Preuve}: Elle est bas\'ee sur la m\'ethode des caract\'eristiques. Nous int\'egrons les \'equations (\ref{equation_n(t,m,a)}) et (\ref{equation_p(t,m,a)}) par rapport \`a l'\^age :
\setlength\arraycolsep{2pt}
\begin{eqnarray}
\frac{\partial }{\partial t} P(t,m)+\frac{\partial }{\partial m} (V(m)P(t,m)) & = & -\gamma (m)P(t,m) +p(t,m,0) \nonumber \\
     &  &  -p(t,m,\overline{\tau}), \nonumber \\
\frac{\partial }{\partial t} N(t,m)+\frac{\partial }{\partial m} (V(m)N(t,m)) & = & -(\delta (m)+\beta (m, N(t,m)))N(t,m) \nonumber \\
     &  &  +n(t,m,0), \nonumber
\end{eqnarray}
en ayant suppos\'e que $\lim_{a \to \infty} n(t,m,a) = 0$, ce qui signifie simplement, d'un point de vue biologique, que toutes les cellules sont condamn\'ees \`a mourir.

En utilisant la m\'ethode des caract\'eristiques nous obtenons
\begin{equation} \label{pcarac}
    p(t,m,a)=\left\{ \begin{array}{ll}
    p(0,\pi_{-t}(m),a-t)\xi (m,t), & \textrm{si} \ 0 \leq t < a, \\
    p(t-a,\pi_{-a}(m),0)\xi (m,a), & \textrm{si} \ t \geq a.
    \end{array} \right.
\end{equation}

Nous pouvons alors \'ecrire, en utilisant (\ref{pcarac}),
\begin{displaymath}
\begin{array}{rcl}
p(t,m,\overline{\tau}) & = & p(0,\pi_{-t}(m),\overline{\tau}-t)\xi (m,t), \\
         & = & \Gamma (\pi_{-t}(m),\overline{\tau} -t) \xi (m,t),
    \end{array}
\end{displaymath}
si $0 \leq t < \overline{\tau}$, et
\begin{displaymath}
\begin{array}{rcl}
p(t,m,\overline{\tau}) & = & p(t-\overline{\tau},\pi_{-\overline{\tau}}(m),0)\xi (m,\overline{\tau}),   \\
         & = & \beta (\pi_{-\overline{\tau}}(m), N(t-\overline{\tau},\pi _{-\overline{\tau}}(m)))N(t-\overline{\tau},\pi _{-\overline{\tau}}(m)) \xi (m,\overline{\tau}),
\end{array}
\end{displaymath}
si $t \geq \overline{\tau}$.

En utilisant la relation (\ref{condbordp_nonlineaire}), nous obtenons bien que $P(t,m)$ v\'erifie les relations (\ref{P(t,m)_1}) et (\ref{P(t,m)_2}).

De m\^eme, d'apr\`es (\ref{condbordn}) et en utilisant (\ref{pcarac}), nous pouvons \'ecrire :
\begin{eqnarray}
n(t,m,0) & = & 2\int_{\underline{\tau}}^{\overline{\tau}} k(m,a)p(t,g^{-1}(m),a)da  \nonumber\\
     & = & 2\int_{\underline{\tau}}^{\overline{\tau}} k(m,a)p(0,\Delta(t,m),a-t)\xi(g^{-1}(m),t)da, \nonumber\\
     & = & 2\xi (g^{-1}(m),t)\int_{\underline{\tau}}^{\overline{\tau}} k(m,a)\Gamma(\Delta(t,m), a-t)da, \nonumber
\end{eqnarray}
si $0 \leq t < \underline{\tau}$, et
\begin{eqnarray}
n(t,m,0) & = & 2\int_{\underline{\tau}}^{t} k(m,a)p(t,g^{-1}(m),a)da \nonumber \\
     &   & +2\int_{t}^{\overline{\tau}} k(m,a)p(t,g^{-1}(m),a)da, \nonumber \\
     & = & 2\int_{\underline{\tau}}^{t} k(m,a)p(t-a,\Delta(a,m),0)\xi (g^{-1}(m),a) da \nonumber \\
     &   & +2\int_{t}^{\overline{\tau}} k(m,a)p(0,\Delta(t,m),a-t)\xi (g^{-1}(m),t)da, \nonumber \\
     & = & 2\int_{\underline{\tau}}^{t} k(m,a)\beta (\Delta(a,m), N(t-a,\Delta(a,m))) N(t-a,\Delta(a,m)) \times \nonumber\\
     &   & \xi (g^{-1} (m),a)da+2\xi (g^{-1}(m),t)\int_{t}^{\overline{\tau}} k(m,a)\Gamma (\Delta(t,m), a-t) da \nonumber,
\end{eqnarray}
si $\underline{\tau} \leq t < \overline{\tau}$, et
\begin{eqnarray}
n(t,m,0) & = & 2\int_{\underline{\tau}}^{\overline{\tau}} k(m,a)p(t-a,\Delta(a,m),0)\xi (g^{-1}(m),a)da, \nonumber \\
     & = & 2\int_{\underline{\tau}}^{\overline{\tau}} k(m,a)\beta (\Delta(a,m), N(t-a,\Delta(a,m))) \times \nonumber \\
     &   & \ \ N(t-a,\Delta(a,m))\xi (g^{-1}(m),a)da. \nonumber
\end{eqnarray}
si $t \geq \overline{\tau}$.

$N(t,m)$ v\'erifie donc bien les \'equations (\ref{N(t,m)_1bis}), (\ref{N(t,m)_1}) et (\ref{N(t,m)_2}).

En remarquant que, pour tout $m\in[g(1),1]$, $\zeta (m,a)=0$ et $k(m,a)=0$, nous obtenons l'\'equation (\ref{N(t,m)_3}).$\Box$

\vspace{2ex}
Nous remarquons tout d'abord que les solutions des \'equations (\ref{N(t,m)_1bis}), (\ref{N(t,m)_1}) et (\ref{P(t,m)_1}) deviennent des conditions initiales pour les \'equations (\ref{N(t,m)_2}) et (\ref{P(t,m)_2}) respectivement. De plus, les solutions de l'\'equation (\ref{N(t,m)_2}) ne d\'ependent pas des solutions de l'\'equation (\ref{P(t,m)_2}). Enfin, si nous connaissons le comportement des solutions de (\ref{N(t,m)_2}) nous pouvons en d\'eduire les solutions de l'\'equation (\ref{P(t,m)_2}) par la m\'ethode des caract\'eristiques. Nous concentrerons donc notre \'etude sur les \'equations d\'ecrivant le comportement des cellules dans la phase de repos.

Nous allons d\'eterminer une formulation faible du probl\`eme (\ref{N(t,m)_1bis})-(\ref{N(t,m)_1})-(\ref{N(t,m)_2}).

\section{Formulation faible} \label{para_formulationfaible}

Soit $\mathcal{C}^{0}[0,g(1)]$ l'espace des fonctions continues sur $[0,g(1)]$, muni de la norme du maximum, et soit $A:\mathcal{D}(A) \subset \mathcal{C}^{0}[0,g(1)] \to \mathcal{C}^{0}[0,g(1)]$ l'op\'erateur d\'efini sur l'ensemble :
\begin{displaymath}
\mathcal{D}(A) = \Big\{u \in \mathcal{C}^{0}[0,g(1)] ; u \in \mathcal{C}^{1}(0,g(1)], \lim_{m \to 0} V(m)u'(m) =0\Big\},
\end{displaymath}
par
\begin{displaymath}
Au(m)=\left\{\begin{array}{ll}
    -(\delta (m) +V'(m))u(m)-V(m)u'(m), & \textrm{si} \ m \in (0,g(1)], \\
    -(\delta (0) +V'(0))u(0), & \textrm{si} \ m =0.
    \end{array} \right.
\end{displaymath}

\vspace{1ex}
\begin{proposition}
$A$ est le g\'en\'erateur infinit\'esimal du $C_{0}$-semi-groupe $(S(t))_{t \geq 0}$ d\'efini, pour $\psi~\in~\mathcal{D}(A)$, $m \in [0,g(1)]$ et $t \geq 0$, par:
\begin{displaymath}
(S(t)\psi )(m)=\psi (h^{-1}(h(m)e^{-t}))exp\bigg\{-\int_{0}^{t} (\delta+V')(h^{-1}(h(m)e^{-s}))ds \bigg\}.
\end{displaymath}
\end{proposition}

\vspace{1ex}
\begin{remarque}
Comme $h^{-1}(h(m)e^{s})=\pi_{s}(m)$ pour tout $m \in [0,1]$ et tout $s \leq 0$, alors, $\forall \psi~\in~\mathcal{D}(A)$, $\forall m \in [0,g(1)]$ et $ \forall t \geq 0$ :
\begin{equation} \label{eq0}
    (S(t)\psi )(m)=\psi (\pi_{-t}(m))K(t,m),
\end{equation}
o\`u
\begin{displaymath}
    K(t,m)=exp\bigg\{-\int_{0}^{t} \delta (\pi_{-s}(m)) +V'(\pi_{-s}(m)) ds \bigg\},
\end{displaymath}
pour tout $m \in [0,g(1)]$ et tout $t \geq 0$.
\end{remarque}
\vspace{1ex}
\textbf{Preuve de la Proposition :} La d\'emonstration de cette proposition n\'ecessite l'utilisation de trois lemmes.

Soit $B:\mathcal{D}(B) \subset \mathcal{C}^{0}[0,g(1)] \to \mathcal{C}^{0}[0,g(1)]$ l'op\'erateur d\'efini par :
\begin{displaymath}
\mathcal{D}(B) = \Big\{u \in \mathcal{C}^{0}[0,g(1)] \, ; u \in \mathcal{C}^{1}(0,g(1)], \lim_{m \to 0} V(m)u'(m) =0\Big\},
\end{displaymath}
\begin{equation} \label{operateurB}
Bu(m)=\left\{\begin{array}{ll}
     -V(m)u'(m), & \textrm{si} \ m \in (0,g(1)], \\
      0,        & \textrm{si} \ m =0.
    \end{array} \right.
\end{equation}

\begin{lemme}
L'op\'erateur $B$ d\'efini par (\ref{operateurB}) sur $\mathcal{D}(B)$ est le g\'en\'erateur infinit\'esimal d'un semi-groupe de contractions  sur $\mathcal{C}^{0}[0,g(1)]$.
\end{lemme}
\vspace{2ex}
\textbf{Preuve :} Nous allons utiliser le th\'eor\`eme de Hille-Yosida. $B$ est un op\'erateur lin\'eaire ferm\'e. De plus, en remarquant que $\mathcal{C}^{1}[0,g(1)] \subset \mathcal{D}(B) \subset \mathcal{C}^{0}[0,g(1)]$ et comme $\mathcal{C}^{1}[0,g(1)]$ est dense dans $\mathcal{C}^{0}[0,g(1)]$, nous obtenons bien que $\mathcal{D}(B)$ est dense dans $\mathcal{C}^{0}[0,g(1)]$. \\
Soient $w \in \mathcal{C}^{0}[0,g(1)]$ et  $\lambda > 0$ donn\'es. Montrons qu'il existe une fonction $u \in \mathcal{D}(B)$ et une seule telle que
\begin{displaymath}
(I-\lambda B)u=w,
\end{displaymath}
o\`u $I$ est l'op\'erateur identit\'e, et que
\begin{displaymath}
|u|_{\infty} \leq |w|_{\infty},
\end{displaymath}
o\`u la norme $|.|_{\infty}$ est d\'efinie par
\begin{displaymath}
|u|_{\infty}=\sup_{0\leq m \leq g(1)}|u(m)|.
\end{displaymath}
Si $u$ existe alors $u$ est solution de l'\'equation diff\'erentielle
\begin{equation} \label{equadiff}
\left\{ \begin{array}{l}
u(m)+\lambda V(m)u'(m)=w(m), \quad \textrm{pour } 0<m\leq g(1), \\
u(0)=w(0).
\end{array} \right.
\end{equation}
Consid\'erons la fonction $h_{\lambda}$ d\'efinie par
\begin{displaymath}
h_{\lambda}(m)=\left\{ \begin{array}{ll}
exp\Big( \displaystyle\int_{g(1)}^{m} \frac{ds}{\lambda V(s)} \Big), & \textrm{pour } m \in (0,g(1)], \\
0,                            & \textrm{pour } m = 0.
\end{array} \right.
\end{displaymath}
Les solutions de l'\'equation (\ref{equadiff}) sont de la forme
\begin{equation} \label{solutionequadiff}
u(m)=\frac{1}{h_{\lambda}(m)}\bigg(k+\int_{g(1)}^{m} \frac{w(s)}{\lambda V(s)} h_{\lambda}(s)ds \bigg),
\end{equation}
pour $m \in (0,g(1)]$ et $k \in \mathbb{R}$.
$h_{\lambda}$ est continue sur $[0,g(1)]$, de classe $\mathcal{C}^{1}$ sur $(0,g(1)]$ et
\begin{displaymath}
h_{\lambda}'(m)=\frac{h_{\lambda}(m)}{\lambda V(m)},
\end{displaymath}
pour $m \in (0,g(1)]$. Donc
\begin{displaymath}
\int_{g(1)}^{m}\frac{h_{\lambda}(s)}{\lambda V(s)}ds = h_{\lambda}(m)-h_{\lambda}(g(1)).
\end{displaymath}
Ainsi, $\lim_{m \to 0}\int_{g(1)}^{m}\frac{h_{\lambda}(s)}{\lambda V(s)}ds =  -h_{\lambda}(g(1)) =  -1$. Par cons\'equent, $\frac{h_{\lambda}(s)}{\lambda V(s)} \in L^{1}(0,g(1))$ et, comme $h_{\lambda}(0)=0$, la fonction $u$ donn\'ee par (\ref{solutionequadiff}) est born\'ee si et seulement si
\begin{displaymath}
k=\int_{0}^{g(1)} \frac{w(s)}{\lambda V(s)}h_{\lambda}(s) ds.
\end{displaymath}
L'unique solution possible de l'\'equation (\ref{equadiff}) est donc
\begin{displaymath}
u(m)=\frac{1}{h_{\lambda}(m)}\bigg(\int_{0}^{m} \frac{w(s)}{\lambda V(s)} h_{\lambda}(s)ds \bigg),
\end{displaymath}
pour $m \in (0,g(1)]$, et $u(0)=w(0)$.

\noindent De plus, si $u \in \mathcal{D}(B)$, alors $|u(m)| \leq
|w|_{\infty}$ pour tout $m \in (0,g(1)]$ et, comme $u(0)=w(0)$,
$|u|_{\infty} \leq |w|_{\infty}$.

\noindent En appliquant le th\'eor\`eme de Hille-Yosida, nous
achevons la d\'emonstration du lemme. $\Box$\bigskip{}

\begin{lemme} \label{generateurB}
Si $\psi \in \mathcal{D}(B)$, alors
\begin{displaymath}
\lim_{t \to 0^{+}} \frac{\psi(h^{-1}(h(m)e^{-t}))-\psi(m)}{t}=-V(m)\psi'(m).
\end{displaymath}
\end{lemme}
\vspace{2ex}
\textbf{Preuve :} Nous effectuons le changement de variables suivant :
\begin{displaymath}
y=h^{-1}(h(m)e^{-t}).
\end{displaymath}
Alors, par d\'efinition de $h$,
\begin{displaymath}
t=\ln \bigg(\frac{h(m)}{h(y)}\bigg)=\int_{y}^{m}\frac{ds}{V(s)}.
\end{displaymath}
Ainsi
\begin{eqnarray}
\lim_{t \to 0^{+}} \frac{\psi(h^{-1}(h(m)e^{-t}))-\psi(m)}{t} & = & \lim_{y \to m} \frac{\psi(y)-\psi(m)}{y-m} \frac{y-m}{\int_{y}^{m}\frac{ds}{V(s)}} \nonumber \\
                                  & = & -V(m)\psi'(m). \nonumber
\end{eqnarray}
Ce qui conclut la d\'emonstration. $\Box$ \bigskip{}

Par d\'efinition du g\'en\'erateur infinit\'esimal d'un semi-groupe, nous d\'edui\-sons de ce lemme que le semi-groupe engendr\'e par $B$, que nous notons $(T(t))_{t \geq 0}$, est donn\'e, pour tout $\psi \in \mathcal{D}(B)$ et tout $t \geq 0$, par
\begin{displaymath}
(T(t)\psi)(m)=\psi(h^{-1}(h(m)e^{-t})), \, \forall m \in [0,g(1)].
\end{displaymath}

\vspace{1ex}
\begin{lemme} \label{generateurB+alphaI}
Si $(T(t))_{t \geq 0}$ est le semi-groupe engendr\'e par $B$ alors le semi-groupe engendr\'e par $B+\alpha I$, avec $\alpha \in \mathcal{C}^{0}[0,g(1)]$, est d\'efini par
\begin{displaymath}
(S(t)\psi )(m)=exp\bigg( \int_{0}^{t} (T(s)\alpha)(m)ds \bigg)(T(t)\psi)(m),
\end{displaymath}
pour tout $\psi \in \mathcal{D}(B)$, tout $m \in [0,g(1)]$ et tout $t \geq 0$.
\end{lemme}
\vspace{2ex}
\textbf{Preuve :} Il suffit de montrer que pour $\psi~\in~\mathcal{D}(B)$,
\begin{displaymath}
\lim_{t \to 0^{+}} \frac{(S(t)\psi )(m)-\psi(m)}{t} = (B+\alpha I)\psi(m).
\end{displaymath}
D'apr\`es la d\'efinition du semi-groupe $(S(t))$ nous avons :
\begin{eqnarray}
\frac{(S(t)\psi )(m)-\psi(m)}{t} & = & \frac{1}{t}\bigg[exp\bigg( \int_{0}^{t} (T(s)\alpha)(m)ds \bigg)(T(t)\psi)(m)-\psi(m)\bigg], \nonumber \\
                 & = & exp\bigg( \int_{0}^{t} (T(s)\alpha)(m)ds \bigg) \bigg(\frac{(T(t)\psi )(m)-\psi(m)}{t}\bigg) \nonumber \\
                 & + & \frac{1}{t}\bigg[exp\bigg( \int_{0}^{t} (T(s)\alpha)(m)ds \bigg)\psi(m)-\psi(m) \bigg]. \nonumber
\end{eqnarray}
En passant \`a la limite sur $t$ nous obtenons le r\'esultat : le
premier terme de droite tend vers $B\psi(m)$ et le deuxi\`eme vers
$(\alpha I)\psi(m)$.  $\Box$\bigskip{}

En utilisant les lemmes \ref{generateurB} et
\ref{generateurB+alphaI}, avec $\alpha = \delta +V'$, nous
obtenons le r\'esultat \'enonc\'e dans la proposition.
$\Box$.\bigskip{}

\begin{proposition}
Soit $\varphi \in\mathcal{C}^{0}([0,\overline{\tau}]\times[0,g(1)])$. La formulation int\'egr\'ee du probl\`eme (\ref{N(t,m)_1bis})-(\ref{N(t,m)_1})-(\ref{N(t,m)_2}) est :
\begin{equation} \label{formulationfaible_nonlin}
\left\{ \begin{array}{rcl}
N(t,m) & = & \varphi(\overline{\tau},\pi_{-(t-\overline{\tau})}(m)) K(t-\overline{\tau},m)+ G(N)(t,m) \\
    &   & \qquad \qquad\qquad \ \ \ -J(N)(t,m) \qquad \qquad si \ t \geq \overline{\tau}, \\
N(t,m) & = & \varphi (t,m), \qquad \qquad \qquad \qquad \qquad \qquad \quad si \ 0 \leq t \leq \overline{\tau},
\end{array} \right.
\end{equation}
pour tout $m \in [0,g(1)]$, avec
\setlength\arraycolsep{2pt}
\begin{eqnarray}
G(N)(t,m) & = & 2\!\!\int_{\overline{\tau}}^{t}\!\bigg(\!\int_{\underline{\tau}}^{\overline{\tau}}\!\! \zeta (\pi_{-(t-s)}(m),a)\beta (\Delta(a,\pi_{-(t-s)}(m)),N(s-a,\Delta(a,\pi_{-(t-s)}(m)))\times \nonumber \\
     &   & \qquad \qquad\qquad N(s-a,\Delta(a,\pi_{-(t-s)}(m))) da \bigg) K(t-s,m) ds, \nonumber
\end{eqnarray}
et
\begin{displaymath}
J(N)(t,m)=\int_{\overline{\tau}}^{t}K(t-s,m)\beta(\pi_{-(t-s)}(m),N(s,\pi_{-(t-s)}(m)))N(s,\pi_{-(t-s)}(m))ds,
\end{displaymath}
pour tous $t \geq \overline{\tau}$ et $m \in [0,g(1)]$.
\end{proposition}
\vspace{2ex}
\textbf{Preuve :} L'\'equation (\ref{N(t,m)_2}) se met sous la forme :
\setlength\arraycolsep{2pt}
\begin{displaymath}
\frac{\partial }{\partial t}N(t,m) = -(\delta (m)+\beta (m,N(t,m))+V'(m))N(t,m)-V(m)\frac{\partial }{\partial m} N(t,m)
\end{displaymath}
\begin{equation} \label{eq1}
+2\int_{\underline{\tau}}^{\overline{\tau}} \zeta (m,a)\beta (\Delta(a,m),N(t-a,\Delta(a,m)))N(t-a,\Delta(a,m))da.
\end{equation}

Soit $u$ l'application d\'efinie pour tout $t \geq 0$ par:
\begin{equation} \label{u(t)}
    u(t)(m)=N(t,m), \qquad \qquad \forall m \in [0,1].
\end{equation}
Comme
\begin{displaymath}
N(t-a,\Delta(a,m))=u(t-a)(\Delta(a,m)):=u_{-a}(t)(\Delta(a,m)),
\end{displaymath}
il vient, d'apr\`es (\ref{eq1}), pour tout $m \in [0,1]$ :
\setlength\arraycolsep{2pt}
\begin{eqnarray}
u'(t)(m) & = & -(\delta (m)+V'(m))u(t)(m) -V(m)(u(t))'(m) -\beta(m,u(t)(m))u(t)(m)\nonumber \\
     &   & +2\int_{\underline{\tau}}^{\overline{\tau}} \zeta (m,a)\beta (\Delta(a,m),u(t-a)(\Delta(a,m)))u_{-a}(t)(\Delta(a,m))da ,\nonumber
\end{eqnarray}
c'est-\`a-dire
\begin{eqnarray*}
u'(t)(m) & = & Au(t)(m)-\beta(m,u(t)(m))u(t)(m) \\
    &   & +2\int_{\underline{\tau}}^{\overline{\tau}} \zeta (m,a)\beta (\Delta(a,m),u(t-a)(\Delta(a,m)))u_{-a}(t)(\Delta(a,m))da.
\end{eqnarray*}
$u(t)$ est donc donn\'e, pour tout $m \in [0,g(1)]$,  par la formule :
\setlength\arraycolsep{2pt}
\begin{eqnarray*}
u(t)(m) & = & S(t-\overline{\tau})\varphi(\overline{\tau})(m)-\displaystyle\int_{\overline{\tau}}^{t} S(t-s)\beta(\pi_{-(t-s)}(m),u(s)(\pi_{-(t-s)}(m)))\times \\
    &   & u(s)(\pi_{-(t-s)}(m))ds + 2\displaystyle\int_{\overline{\tau}}^{t} S(t-s)\Big(\displaystyle\int_{\underline{\tau}}^{\overline{\tau}} \zeta (m,a) \times \\
    &   & \beta (\Delta(a,m),u(t-a)(\Delta(a,m)))u_{-a}(s)(\Delta(a,m))da\Big) ds,
\end{eqnarray*}
si $t \geq \overline{\tau}$, et $u(t)(m)=\varphi(t)(m) := \varphi(t,m)$ si $t \leq \overline{\tau}$, o\`u $\varphi(t) \in \mathcal{C}^{0}([0,g(1)])$ pour tout $t \geq 0$.

\vspace{1ex}
D'apr\`es (\ref{eq0}), si $t \geq \overline{\tau}$,
\setlength\arraycolsep{2pt}
\begin{eqnarray}
u(t)(m) & = & \varphi(\overline{\tau},\pi_{-(t-\overline{\tau})}(m)) K(t-\overline{\tau},m)-\int_{\overline{\tau}}^{t} \beta(\pi_{-(t-s)}(m),u(s)(\pi_{-(t-s)}(m)))\times  \nonumber  \\
        &  & u(s)(\pi_{-(t-s)}(m))K(t-s,m)ds + 2\int_{\overline{\tau}}^{t}\bigg(\int_{\underline{\tau}}^{\overline{\tau}} \zeta(\pi_{-(t-s)}(m),a) \times \nonumber \\
        &   &  \beta (\Delta(a,m),u(t-a)(\Delta(a,m)))u_{-a}(s)(\Delta(a,\pi_{-(t-s)}(m))) da \bigg)K(t-s,m) ds. \nonumber
\end{eqnarray}
En utilisant (\ref{u(t)}) nous obtenons alors (\ref{formulationfaible_nonlin}). $\Box$

\section{Existence et unicit\'e} \label{para_existenceunicite}

Nous allons \`a pr\'esent \'etablir l'existence et l'unicit\'e des solutions du probl\`eme (\ref{formulationfaible_nonlin}), en utilisant une m\'ethode it\'erative et la m\'ethode des pas, avant de montrer que l'unicit\'e des solutions ne d\'epend en fait que des cellules de petites maturit\'es.

Nous faisons l'hypoth\`ese suivante :

\textbf{(H$_{1}$)} Il existe $l\in\mathbb{R}^{+}$ tel que $\forall m \in [0,g(1)], \forall x_{1},x_{2} \in \mathbb{R}^{+},$
\begin{displaymath}
|x_{1}\beta(m,x_{1})-x_{2}\beta(m,x_{2})| \leq l|x_{1}-x_{2}|.
\end{displaymath}

\begin{theoreme} \label{existenceunicite_nonlin}
Sous l'hypoth\`ese \emph{\textbf{(H$_{1}$)}}, pour toute donn\'ee initiale $\varphi \in \mathcal{C}^{0}([0,\overline{\tau}]\times[0,g(1)])$, le probl\`eme (\ref{formulationfaible_nonlin}) poss\`ede une unique solution qui est continue.
\end{theoreme}
\vspace{2ex}
\textbf{Preuve :} Soit $\varphi \in \mathcal{C}^{0}([0,\overline{\tau}]\times[0,g(1)])$ donn\'ee. Si $N$ est solution du probl\`eme (\ref{formulationfaible_nonlin}) alors
\begin{equation}
N(t,m)=\varphi(\overline{\tau},\pi_{-(t-\overline{\tau})}(m)) K(t-\overline{\tau},m) + G(N)(t,m) - J(N)(t,m),
\end{equation}
pour tout $t \geq \overline{\tau}$ et tout $m \in [0,g(1)]$, et $N(t,m)=\varphi(t,m)$ pour $t \in [0,\overline{\tau}]$ et $m \in [0,g(1)]$.

Soit $t \in [0,\overline{\tau}+\underline{\tau}]$. Nous d\'efinissons la suite $(N_{n})_{n \in \mathbb{N}}$ de la fa\c con suivante : \\
si $t \in [0,\overline{\tau}]$ alors $N_{n}(t,m)=\varphi(t,m)$, $\forall m \in [0,g(1)]$, et si $t \in [\overline{\tau},\overline{\tau}+\underline{\tau}]$ alors, $\forall m \in [0,g(1)]$,
\begin{displaymath}
\left\{ \begin{array}{rcl}
N_{0}(t,m) & = & \varphi(\overline{\tau},\pi_{-(t-\overline{\tau})}(m)) K(t-\overline{\tau},m) + G(\varphi)(t,m), \\
N_{n}(t,m) & = & N_{0}(t,m) - J(N_{n-1})(t,m).
\end{array} \right.
\end{displaymath}
Nous avons alors:
\begin{displaymath}
|N_{1}(t,m)-N_{0}(t,m)|=|J(N_{0})(t,m)|.
\end{displaymath}
Les fonctions $\delta$ et $V'$ \'etant continues sur $[0,g(1)]$, l'application $K$ est continue sur $[0,+\infty)\times[0,g(1)]$ et donc born\'ee sur $[\overline{\tau},\overline{\tau}+\underline{\tau}]\times[0,g(1)]$. Notons
\begin{displaymath}
\alpha = \max_{(t,m) \in [\overline{\tau},\overline{\tau}+\underline{\tau}]\times[0,g(1)]} |K(t,m)|.
\end{displaymath}
Alors,
\begin{eqnarray}
|N_{1}(t,m)-N_{0}(t,m)| & \leq & \alpha \int_{\overline{\tau}}^{t} |\beta(\pi_{-(t-s)}(m),N_{0}(s,\pi_{-(t-s)}(m)))|\times \nonumber \\
            &      & \qquad \qquad |N_{0}(s,\pi_{-(t-s)}(m))|ds, \nonumber \\
            & \leq & \alpha l \int_{\overline{\tau}}^{t} |N_{0}(s,\pi_{-(t-s)}(m))| ds. \nonumber
\end{eqnarray}
$\varphi$, $K$ et $G(\varphi)$ \'etant continues, $N_{0}$ est continue sur $[\overline{\tau},\overline{\tau}+\underline{\tau}]\times[0,g(1)]$ :
\begin{displaymath}
\exists M \geq 0,\ \forall m \in [0,g(1)],\ \forall t \in [\overline{\tau},\overline{\tau}+\underline{\tau}],\ |N_{0}(t,m)| \leq M.
\end{displaymath}
Alors,
\begin{displaymath}
|N_{1}(t,m)-N_{0}(t,m)| \leq M\alpha l (t-\overline{\tau}).
\end{displaymath}
Montrons par r\'ecurrence le r\'esultat suivant :
\begin{displaymath}
(\textrm{\textbf{H}}_{n}) \ \forall t \in [\overline{\tau},\overline{\tau}+\underline{\tau}], \forall m \in [0,g(1)], |N_{n}(t,m)-N_{n-1}(t,m)| \leq M(\alpha l)^{n} \frac{(t-\overline{\tau})^{n}}{n!}.
\end{displaymath}
La relation est vraie pour $n=0$, nous venons de le v\'erifier. Supposons $(\textrm{\textbf{H}}_{n})$ satisfaite et montrons que $(\textrm{\textbf{H}}_{n+1})$ est r\'ealis\'ee. \\
Soient $t \in [\overline{\tau},\overline{\tau}+\underline{\tau}]$ et $m \in [0,g(1)]$. Alors,
\begin{displaymath}
|N_{n+1}(t,m)-N_{n}(t,m)|=|J(N_{n})(t,m)-J(N_{n-1})(t,m)|,
\end{displaymath}
donc, en utilisant l'hypoth\`ese \textbf{(H$_{1}$)},
\begin{eqnarray*}
|N_{n+1}(t,m)-N_{n}(t,m)| & \leq & \alpha l \int_{\overline{\tau}}^{t}|N_{n}(s,\pi_{-(t-s)}(m)) \\
            &   & \qquad \qquad  -  N_{n-1}(s,\pi_{-(t-s)}(m))| ds.
\end{eqnarray*}
En utilisant l'hypoth\`ese de r\'ecurrence, nous obtenons
\begin{eqnarray}
|N_{n+1}(t,m)-N_{n}(t,m)| & \leq & M(\alpha l)^{n+1} \int_{\overline{\tau}}^{t} \frac{(s-\overline{\tau})^{n}}{n!} ds, \nonumber \\
            & \leq & M(\alpha l)^{n+1} \frac{(t-\overline{\tau})^{n+1}}{(n+1)!}. \nonumber
\end{eqnarray}
Nous obtenons bien le r\'esultat escompt\'e.\\
Par cons\'equent, $N:=\lim_{n \to \infty} N_{n}$ existe uniform\'ement sur $[\overline{\tau},\overline{\tau}+\underline{\tau}]$ et est continue.

Montrons que $N$ est solution de (\ref{formulationfaible_nonlin}) pour $t \in [\overline{\tau},\overline{\tau}+\underline{\tau}]$ et $m \in [0,g(1)]$, c'est-\`a-dire
\begin{displaymath}
N(t,m)=N_{0}(t,m) - J(N)(t,m),
\end{displaymath}
car $G(N)(t,m)=G(\varphi)(t,m)$, pour tout $t \in [\overline{\tau},\overline{\tau}+\underline{\tau}]$ et tout $m \in [0,g(1)]$. En effet, si $a \in [\underline{\tau},\overline{\tau}]$ et $s \in [\overline{\tau},t]$ alors $s-a \in [0,t-\underline{\tau}] \subset [0,\overline{\tau}]$ et donc
\begin{displaymath}
N(s-a,\Delta(a,m))=\varphi(s-a,\Delta(a,m)),
\end{displaymath}
pour tout $m \in [0,g(1)]$.
Posons
\begin{displaymath}
v(t,m)=|N(t,m)-N_{0}(t,m)-J(N)(t,m)|,
\end{displaymath}
pour tout $t \in [\overline{\tau},\overline{\tau}+\underline{\tau}]$ et tout $m \in [0,g(1)]$ et montrons que $v \equiv 0$.\\
En remarquant que
\begin{displaymath}
N(t,m)-N_{n}(t,m)=\sum_{p=n}^{\infty} (N_{p+1}(t,m)-N_{p}(t,m)),
\end{displaymath}
alors
\begin{eqnarray}
v(t,m)  & \leq & |N(t,m)-N_{n+1}(t,m)|+|J(N_{n})(t,m)-J(N)(t,m)|, \nonumber \\
    & \leq & |N(t,m)-N_{n+1}(t,m)| \nonumber \\
    &      & \quad +\alpha l \int_{\overline{\tau}}^{t} |N_{n}(s,\pi_{-(t-s)}(m))-N(s,\pi_{-(t-s)}(m))| ds, \nonumber \\
    & \leq & \sum_{p=n+1}^{\infty} |N_{p+1}(t,m)-N_{p}(t,m)| \nonumber \\
    &      & \ +\alpha l \int_{\overline{\tau}}^{t} \sum_{p=n}^{\infty} |N_{p+1}(s,\pi_{-(t-s)}(m))-N_{p}(s,\pi_{-(t-s)}(m))| ds, \nonumber \\
    & \leq & M \sum_{p=n+1}^{\infty} (\alpha l)^{p+1}\frac{(t-\overline{\tau})^{p+1}}{(p+1)!} + M\alpha l \sum_{p=n}^{\infty} \int_{\overline{\tau}}^{t} (\alpha l)^{p+1} \frac{(s-\overline{\tau})^{p+1}}{(p+1)!}ds, \nonumber \\
    & \leq & M \sum_{p=n+1}^{\infty} (\alpha l)^{p+1}\frac{(t-\overline{\tau})^{p+1}}{(p+1)!} + M \sum_{p=n}^{\infty} (\alpha l)^{p+2} \frac{(t-\overline{\tau})^{p+2}}{(p+2)!}, \nonumber \\
    & \leq & 2M \sum_{p=n+1}^{\infty} (\alpha l)^{p+1}\frac{(t-\overline{\tau})^{p+1}}{(p+1)!}. \nonumber
\end{eqnarray}
En faisant tendre $n$ vers $+\infty$, nous obtenons :
\begin{displaymath}
v(t,m)=0, \qquad \forall t \in [\overline{\tau},\overline{\tau}+\underline{\tau}] \textrm{ et } \forall m \in [0,g(1)].
\end{displaymath}
Ainsi, $N(t,m)$ est solution de (\ref{formulationfaible_nonlin}) pour $t \in [\overline{\tau},\overline{\tau}+\underline{\tau}]$ et $ m \in [0,g(1)]$.

Montrons que $N$ est unique.\\
Soit $Z(t,m)$  une autre solution de (\ref{formulationfaible_nonlin}) sur $[\overline{\tau},\overline{\tau}+\underline{\tau}]\times[0,g(1)]$ :
\begin{displaymath}
Z(t,m)=N_{0}(t,m) - J(Z)(t,m),
\end{displaymath}
pour tout $t \in [\overline{\tau},\overline{\tau}+\underline{\tau}]$ et tout $ m \in [0,g(1)]$.\\
Alors, $\forall t \in [\overline{\tau},\overline{\tau}+\underline{\tau}]$ et $\forall m \in [0,g(1)]$,
{\setlength\arraycolsep{2pt}
\begin{eqnarray*}
|N(t,m)-Z(t,m)| & = & |J(N)(t,m)-J(Z)(t,m)|,  \\
        & \leq & \alpha l \int_{\overline{\tau}}^{t} |N(s,\pi_{-(t-s)}(m))-Z(s,\pi_{-(t-s)}(m))| ds.
\end{eqnarray*}
Pour $t \in [\overline{\tau},\overline{\tau}+\underline{\tau}]$ et $m \in [0,g(1)]$ nous posons :
\begin{displaymath}
w(\sigma)=|N(\sigma,\pi_{-(t-\sigma)}(m))-Z(\sigma,\pi_{-(t-\sigma)}(m))|,
\end{displaymath}
pour tout $\sigma \in [\overline{\tau},t]$.\\
Alors,
{\setlength\arraycolsep{2pt}
\begin{eqnarray*}
w(\sigma) & = & |J(N)(\sigma,\pi_{-(t-\sigma)}(m))-J(Z)(\sigma,\pi_{-(t-\sigma)}(m))|,  \\
      & \leq & \alpha l \int_{\overline{\tau}}^{\sigma} |N(s,\pi_{-(\sigma-s)} \circ \pi_{-(t-\sigma)}(m))  \\
      &      & \qquad \quad -Z(s,\pi_{-(\sigma-s)} \circ \pi_{-(t-\sigma)}(m))| dr.
\end{eqnarray*}
Comme
\begin{displaymath}
\pi_{-(\sigma-s)} \circ \pi_{-(t-\sigma)}=\pi_{-(t-s)},
\end{displaymath}
alors
\begin{eqnarray}
w(\sigma) & \leq & \alpha l \int_{\overline{\tau}}^{\sigma} |N(s,\pi_{-(t-s)}(m))-Z(s,\pi_{-(t-s)}(m))| dr, \nonumber \\
      & \leq & \alpha l \int_{\overline{\tau}}^{\sigma} w(s) ds. \nonumber
\end{eqnarray}
D'apr\`es le lemme de Gronwall, $\forall \sigma \in [\overline{\tau},t]$, $w(\sigma)=0$, donc $\forall \sigma \in [\overline{\tau},\overline{\tau}+\underline{\tau}]$, $w(\sigma)=0$. Ainsi, $\forall \sigma \in [\overline{\tau},\overline{\tau}+\underline{\tau}]$, $\forall m \in [0,g(1)]$,
\begin{displaymath}
N(\sigma,\pi_{-(t-\sigma)}(m))=Z(\sigma,\pi_{-(t-\sigma)}(m)).
\end{displaymath}
Donc $N(t,m)=Z(t,m)$, $\forall t \in [\overline{\tau},\overline{\tau}+\underline{\tau}]$ et $\forall m \in [0,g(1)]$.\\
Le probl\`eme (\ref{formulationfaible_nonlin}) poss\`ede donc une unique solution sur $[0,\overline{\tau}+\underline{\tau}]\times[0,g(1)]$ qui est continue. Notons $N^{*}$ cette solution.

En utilisant la suite $(N_{n})_{n \in \mathbb{N}}$ d\'efinie de la fa\c con suivante : \\
si $t \in [0,\overline{\tau}+\underline{\tau}]$ alors $N_{n}(t,m)=N^{*}(t,m)$, $\forall m \in [0,g(1)]$, et si $t \in [\overline{\tau}+\underline{\tau},\overline{\tau}+2\underline{\tau}]$ alors, $\forall m \in [0,g(1)]$,
\begin{displaymath}
\left\{ \begin{array}{rcl}
N_{0}(t,m) & = & \varphi(\overline{\tau},\pi_{-(t-\overline{\tau})}(m)) K(t-\overline{\tau},m) + G(N^{*})(t,m), \\
N_{n}(t,m) & = & N_{0}(t,m) - J(N_{n-1})(t,m),
\end{array} \right.
\end{displaymath}
nous montrons de la m\^eme mani\`ere que le probl\`eme (\ref{formulationfaible_nonlin}) poss\`ede une unique solution, continue, sur $[0,\overline{\tau}+2\underline{\tau}]\times[0,g(1)]$. Par la m\'ethode des pas, l'existence et l'unicit\'e d'une solution continue sur $[0,+\infty)\times[0,g(1)]$ sont \'etablies. $\Box$ \\
\vspace{1ex}

Avant d\'etablir le r\'esultat d'unicit\'e d\'ependant uniquement des cellules de petites maturit\'es, nous montrons le r\'esultat de positivit\'e suivant.

\begin{proposition}
Nous supposons que $\delta(m)+V'(m)>0$, pour tout $m\in[0,g(1)]$. Si, $\forall (t,m)\in[0,\overline{\tau}]\times[0,g(1)]$, $N(t,m)\geq 0$, alors $N(t,m)\geq 0$, pour tout $t\geq\overline{\tau}$ et tout $m\in[0,g(1)]$.
\end{proposition}
\vspace{1ex}
\textbf{Preuve :} Nous donnons une id\'ee de la d\'emonstration. Nous utilisons la m\'ethode des caract\'eristiques et la m\'ethode des pas.\\
Nous notons $N=\varphi$ sur $[0,\overline{\tau}]\times[0,g(1)]$ et nous supposons $\varphi\geq 0$. Supposons alors qu'il existe $(t_{0},m_{0})\in[\overline{\tau},\overline{\tau}+ \underline{\tau}]\times[0,g(1)]$ tel que $N(t_{0},m_{0})<0$. Consid\'erons alors l'\'equation diff\'erentielle
\begin{displaymath}
\left\{ \begin{array}{rcll}
m'(t)&=&V(m(t)),& \ t\in[\overline{\tau},\overline{\tau}+ \underline{\tau}],\\
m(t_{0})&=&m_{0},
\end{array} \right.
\end{displaymath}
et notons
\begin{displaymath}
u(t)=N(t,m(t)), \ t\in[\overline{\tau},\overline{\tau}+ \underline{\tau}].
\end{displaymath}
Alors
\begin{displaymath}
u'(t)= -(\delta (m(t))+V'(m(t))+\beta (m(t),u(t)))u(t)
\end{displaymath}
\begin{displaymath}
+2\int_{\underline{\tau}}^{\overline{\tau}} \zeta (m(t),a)\beta (\Delta(a,m(t)),\varphi(t-a,\Delta(a,m(t))))\varphi(t-a,\Delta(a,m(t)))da.
\end{displaymath}
Donc $u'(t_{0})>0$. En effet, $u(t_{0})<0$ et
\begin{displaymath}
\int_{\underline{\tau}}^{\overline{\tau}} \zeta (m(t_{0}),a)\beta (\Delta(a,m(t_{0})),\varphi(t_{0}-a,\Delta(a,m(t_{0}))))\varphi(t_{0}-a,\Delta(a,m(t_{0})))da\geq 0.
\end{displaymath}
Nous aboutissons \`a une contradiction. Donc $N(t,m)\geq 0$ sur $[\overline{\tau},\overline{\tau}+ \underline{\tau}]\times[0,g(1)]$. Par la m\'ethode des pas, nous obtenons ce r\'esultat pour $t\geq\overline{\tau}$.$\Box$

\vspace{2ex}
Nous allons \`a pr\'esent \'etablir un r\'esultat qui souligne l'influence des cellules de petites maturit\'es (cellules souches) sur l'\'evolution de la population totale de cellules.

Nous supposerons, dans la suite, qu'il existe $C\geq 0$ tel que pour tout $m\in(0,g(1)]$
\begin{equation} \label{integralefinie}
 \int_{m}^{g^{-1}(m)} \frac{ds}{V(s)} \leq C.
\end{equation}

\begin{exemple}
Si $V(s)\underset{s\to 0}{\sim}\alpha s$, avec $\alpha>0$, alors
\begin{displaymath}
\int_{m}^{g^{-1}(m)} \frac{ds}{V(s)} \leq C \ \textrm{ si et seulement si } \ g'(0) > 0.
\end{displaymath}
\end{exemple}

\vspace{1ex}
\begin{theoreme} \label{unicite_nonlin}
Nous supposons l'hypoth\`ese (\ref{integralefinie}) v\'erifi\'ee.\\
Soit $\tau_{0} =\sup_{m>0} \Big( \int_{m}^{g^{-1}(m)} \frac{ds}{V(s)} \Big)$. Supposons que $N_{1}(t,m)$ et $N_{2}(t,m)$ sont solutions de (\ref{formulationfaible_nonlin}) pour des donn\'ees initiales $\varphi_{1}$ et $\varphi_{2}$ respectivement. Supposons que $\underline{\tau} > \tau_{0}$ et qu'il existe  $b \in (0,h^{-1}(e^{-\underline{\tau}}))$ tel que
\begin{displaymath}
\varphi_{1}(t,m)=\varphi_{2}(t,m), \qquad \textrm{pour } m \in [0,b] \textrm{ et } t \in [0,\overline{\tau}].
\end{displaymath}
Alors, il existe $\overline{t} \geq 0$ tel que
\begin{displaymath}
N_{1}(t,m)=N_{2}(t,m), \qquad \textrm{pour } m \in [0,g(1)] \textrm{ et } t \geq \overline{t}.
\end{displaymath}
\end{theoreme}
\vspace{2ex}
\textbf{Preuve :} Remarquons que la condition $\underline{\tau} > \tau_{0}$ implique, d'apr\`es le lemme \ref{fonctiondelta}, que $\Delta(a,m) < m$ pour tout $m \in [0,g(1)]$ et tout $a \in [\underline{\tau},\overline{\tau}]$, et que $h^{-1}(e^{-\underline{\tau}})<g(1)$. De plus, comme $\underline{\tau} > 0$, $h^{-1}(e^{-\underline{\tau}}) > 0$.

Montrons tout d'abord que $N_{1}(t,m)=N_{2}(t,m)$, pour $ m \in [0,b]$ et $ t \geq 0$.\\
Soit $ m \in [0,b]$. Si $t \in [0, \overline{\tau}]$, le r\'esultat est \'evident.\\
Soit $ t \geq \overline{\tau}$. Rappelons que si $N$ est solution du probl\`eme (\ref{formulationfaible_nonlin}) alors
\begin{displaymath}
N(t,m)=\varphi(\overline{\tau},\pi_{-(t-\overline{\tau})}(m)) K(t-\overline{\tau},m) + G(N)(t,m) - J(N)(t,m),
\end{displaymath}
pour tout $t \geq \overline{\tau}$ et tout $m \in [0,g(1)]$.\\
Comme $t-\overline{\tau} \geq 0$,
\begin{displaymath}
\pi_{-(t-\overline{\tau})}(m) < m \leq b.
\end{displaymath}
Donc
\begin{displaymath}
\varphi_{1}(\overline{\tau},\pi_{-(t-\overline{\tau})}(m))=\varphi_{2}(\overline{\tau},\pi_{-(t-\overline{\tau})}(m)).
\end{displaymath}
Soit $t \in [\overline{\tau},\overline{\tau}+\underline{\tau}]$. Alors, $s-a \in [0,t-\underline{\tau}] \subset [0,\overline{\tau}]$ et
\begin{displaymath}
\Delta(a,\pi_{-(t-s)}(m)) < \pi_{-(t-s)}(m) < m \leq b.
\end{displaymath}
Donc $G(N_{1})(t,m)=G(N_{2})(t,m)$. Par cons\'equent,
\begin{eqnarray}
|N_{1}(t,m)-N_{2}(t,m)| & = & |J(N_{1})(t,m)-J(N_{2})(t,m)| \nonumber \\
            & \leq & \alpha \tilde{l} \int_{\overline{\tau}}^{t} |N_{1}(s,\pi_{-(t-s)}(m))-N_{2}(s,\pi_{-(t-s)}(m))| ds. \nonumber
\end{eqnarray}
Par un raisonnement analogue \`a celui pr\'esent\'e dans la d\'emonstration pr\'ec\'e\-dente et en utilisant le lemme de Gronwall nous obtenons
\begin{displaymath}
N_{1}(t,m)=N_{2}(t,m),
\end{displaymath}
pour tout $t \in [\overline{\tau},\overline{\tau}+\underline{\tau}]$ et tout $m \in [0,b]$. Par la m\'ethode des pas, nous obtenons le r\'esultat pour $t \geq 0$.

Soit $\alpha$ l'application d\'efinie de $[0,g(1)]$ dans $[0,g(1)]$ par
\begin{displaymath}
\alpha(m)=\Delta(\underline{\tau},m).
\end{displaymath}
D'apr\`es le lemme \ref{fonctiondelta}, $\alpha$ est une application croissante, atteignant son maximum en $g(1)$ et v\'erifiant $\alpha(m) < m$ pour tout $m \in (0,g(1)]$. $\alpha$ est donc inversible et nous notons $\Lambda$ son inverse. Comme $\alpha(g(1))=h^{-1}(e^{-\underline{\tau}})$, $\Lambda$ est d\'efinie sur $[0,h^{-1}(e^{-\underline{\tau}})]$. Nous pouvons prolonger $\Lambda$ par continuit\'e sur $[h^{-1}(e^{-\underline{\tau}}),g(1)]$ en posant
\begin{displaymath}
\Lambda(m)=g(1), \textrm{ pour tout } m \in [h^{-1}(e^{-\underline{\tau}}),g(1)].
\end{displaymath}
$\alpha$ \'etant strictement croissante et positive sur $[0,g(1)]$, $\Lambda$ est strictement croissante sur $[0,h^{-1}(e^{-\underline{\tau}})]$.\\
Consid\'erons la suite $(b_{n})_{n \in \mathbb{N}}$ d\'efinie par
\begin{displaymath}
\left\{ \begin{array}{rcll}
b_{n+1} & = & \Lambda(b_{n}), & \textrm{pour } n \in \mathbb{N}, \\
b_{0}   & = & b.          &
\end{array} \right.
\end{displaymath}
La suite $(b_{n})_{n \in \mathbb{N}}$ est croissante ; en effet, $\forall n \in \mathbb{N}$,
\begin{displaymath}
b_{n}=\alpha(b_{n+1})=\Delta(\underline{\tau},b_{n+1}) < b_{n+1}.
\end{displaymath}
Nous d\'efinissons aussi la suite $(t_{n})_{n \in \mathbb{N}}$ par
\begin{displaymath}
t_{n}=\ln\bigg(\frac{h(b_{n})}{h(b)}\bigg) + n\overline{\tau},
\end{displaymath}
c'est-\`a-dire
\begin{displaymath}
\left\{ \begin{array}{rcll}
t_{n+1} & = & t_{n}+\overline{\tau}+\ln\Big(\frac{h(b_{n+1})}{h(b_{n})}\Big), & \textrm{pour } n \in \mathbb{N}, \\
t_{0}   & = & 0.          &
\end{array} \right.
\end{displaymath}
La suite $(t_{n})_{n \in \mathbb{N}}$ est croissante.

Montrons par r\'ecurrence que
\begin{displaymath}
(H_{n}) \quad \textrm{Si } m \in [0,b_{n}] \textrm{ et } t\geq t_{n} \ \textrm{ alors } \ N_{1}(t,m)=N_{2}(t,m).
\end{displaymath}
Le r\'esultat est vrai pour $n=0$, nous venons de le montrer. \\
Supposons $(H_{n})$ v\'erifi\'ee. Soient $m \in [0,b_{n+1}]$ et $t\geq t_{n+1}$. Si $t\geq t_{n+1}$ alors $t\geq t_{n}+\overline{\tau}$. Nous pouvons r\'e\'ecrire les solutions $N_{i}$, $i=1,2$, du probl\`eme (\ref{formulationfaible_nonlin}) de la fa\c con suivante :
\setlength\arraycolsep{2pt}
\begin{displaymath}
\begin{array}{rcl}
N_{i}(t,m) & = & N_{i}(t_{n}+\overline{\tau},\pi_{-(t-t_{n}-\overline{\tau})}(m)) K(t-t_{n}-\overline{\tau},m)-J(N_{i})(t,m) \\
       & + & 2\displaystyle\int_{t_{n}+\overline{\tau}}^{t}\Big(\int_{\underline{\tau}}^{\overline{\tau}} \zeta(\pi_{-(t-s)}(m),a)N_{i}(s-a,\Delta(a,\pi_{-(t-s)}(m))) \times   \\
       &   & \beta(\Delta(a,\pi_{-(t-s)}(m)),N_{i}(s-a,\Delta(a,\pi_{-(t-s)}(m))))da \Big) K(t-s,m) ds,  \\
\end{array}
\end{displaymath}
pour $t \geq t_{n}+\overline{\tau}$.\\
Montrons que $N_{1}(t_{n}+\overline{\tau},\pi_{-(t-t_{n}-\overline{\tau})}(m))=N_{2}(t_{n}+\overline{\tau},\pi_{-(t-t_{n}-\overline{\tau})}(m))$.\\
$t\geq t_{n+1}$ donc
\begin{displaymath}
t-t_{n}-\overline{\tau} \geq \ln\bigg(\frac{h(b_{n+1})}{h(b_{n})}\bigg),
\end{displaymath}
ou encore,
\begin{displaymath}
e^{-(t-t_{n}-\overline{\tau})} \leq \frac{h(b_{n})}{h(b_{n+1})}.
\end{displaymath}
D'o\`u
\begin{displaymath}
h(m)e^{-(t-t_{n}-\overline{\tau})} \leq h(m)\frac{h(b_{n})}{h(b_{n+1})}\leq h(b_{n+1})\frac{h(b_{n})}{h(b_{n+1})}=h(b_{n}),
\end{displaymath}
car $h$ est croissante. Ainsi,
\begin{displaymath}
\pi_{-(t-t_{n}-\overline{\tau})}(m)=h^{-1}(h(m)e^{-(t-t_{n}-\overline{\tau})}) \leq b_{n}.
\end{displaymath}
En utilisant l'hypoth\`ese de r\'ecurrence, nous obtenons le r\'esultat souhait\'e.\\
Enfin, comme $s-a \geq t_{n}$ et
\begin{displaymath}
\Delta(a,\pi_{-(t-s)}(m)) \leq \Delta(a,m) \leq \Delta(\underline{\tau},m) \leq \Delta(\underline{\tau},b_{n+1})= b_{n},
\end{displaymath}
nous obtenons
\begin{displaymath}
|N_{1}(t,m)-N_{2}(t,m)| = |J(N_{1})(t,m)-J(N_{2})(t,m)|.
\end{displaymath}
Par un raisonnement analogue au pr\'ec\'edent, utilisant le lemme de Gronwall, nous v\'erifions l'hypoth\`ese $(H_{n+1})$.

Remarquons que, par d\'efinition de $\Lambda$, il existe $M \in \mathbb{N^{*}}$ tel que
\begin{displaymath}
b_{M} < g(1)=b_{M+1}.
\end{displaymath}
Alors, pour tout $m \in [0,b_{M+1}]=[0,g(1)]$ et tout $t\geq t_{M+1}$,
\begin{displaymath}
N_{1}(t,m)=N_{2}(t,m).
\end{displaymath}
Nous pouvons poser $\overline{t}=t_{M+1}$ et le th\'eor\`eme est alors d\'emontr\'e.\\
$\Box$

\vspace{2ex}
Nous d\'eduisons de ce th\'eor\`eme le corollaire suivant :

\begin{corollaire}
Nous supposons l'hypoth\`ese (\ref{integralefinie}) v\'erifi\'ee.\\
Soit $\tau_{0} =\sup_{m>0} \Big( \int_{m}^{g^{-1}(m)} \frac{ds}{V(s)} \Big)$. Soit $\varphi\in\mathcal{C}^{0}([0,\overline{\tau}]\times[0,g(1)])$ et $N^{\varphi}$ la solution du probl\`eme (\ref{formulationfaible_nonlin}) associ\'ee \`a la donn\'ee initiale $\varphi$. Supposons que $\underline{\tau} > \tau_{0}$ et qu'il existe  $b \in (0,h^{-1}(e^{-\underline{\tau}}))$ tel que $\varphi(t,m)=0$, pour $m \in [0,b]$ et $t \in [0,\overline{\tau}]$. Alors, il existe $\overline{t} \geq 0$ tel que $N^{\varphi}(t,m)=0$, pour $m \in [0,g(1)]$ et $t \geq \overline{t}$.
\end{corollaire}

\vspace{2ex} Ce r\'esultat met en avant l'influence des cellules
de petites maturit\'es (cellules souches) sur la population future
de cellules. Nous remarquons que dans le cas d'une absence de
cellules souches dans le milieu initial (ce qui est le cas d'une
an\'emie aplasique), l'extinction de la population apr\`es un
temps fini $\overline{t}$ est in\'evitable, les cellules ne
parviennent pas \`a survivre.

\section{Invariance -- Application du th\'eor\`eme \ref{unicite_nonlin}} \label{sectioninvariance}

Dans toute cette partie, nous notons $N^{\varphi}$ la solution du probl\`eme (\ref{formulationfaible_nonlin}) pour une donn\'ee initiale $\varphi\in\mathcal{C}^{0}([0,\overline{\tau}]\times[0,g(1)])$. Pour tout r\'eel $b\in(0,g(1)]$, nous d\'efinissons $||.||_{b}$ par
\begin{displaymath}
\forall u\in \mathcal{C}^{0}([0,\overline{\tau}]\times[0,g(1)]), \ \ ||u||_{b}:=\sup \{ |u(t,m)| ; (t,m)\in[0,\overline{\tau}]\times[0,b] \}.
\end{displaymath}

Nous supposons que
\begin{displaymath}
I:= \inf_{m\in[0,g(1)]} (\delta(m)+V'(m)) > 0,
\end{displaymath}
et nous notons
\begin{displaymath}
\tilde{\zeta}:= \sup \{ |\zeta(m,a)| ; (m,a)\in[0,g(1)]\times[0,\overline{\tau}] \}.
\end{displaymath}

\begin{proposition} \label{propositioninvariance}
Supposons que $l(2(\overline{\tau}-\underline{\tau})\tilde{\zeta}+1)<I$. Alors, pour tout $ b\in(0,g(1)]$ et toute fonction $\varphi\in\mathcal{C}^{0}([0,\overline{\tau}]\times[0,g(1)])$,
\begin{displaymath}
\forall t\geq \overline{\tau}, \ \ \forall m\in[0,b], \quad |N^{\varphi}(t,m)|\leq ||\varphi||_{b}.
\end{displaymath}
\end{proposition}
\vspace{1ex}
\textbf{Preuve :} Soient $b\in(0,g(1)]$ et $\varphi\in\mathcal{C}^{0}([0,\overline{\tau}]\times[0,g(1)])$. Pour $t\geq 0$ et $m\in[0,b]$, nous d\'efinissons la suite $(N_{n})_{n\in\mathbb{N}}$ de la fa\c con suivante : $\forall n\in\mathbb{N}$, $N_{n}(t,m)=\varphi(t,m)$ si $t\in[0,\overline{\tau}]$, et, si $t\geq\overline{\tau}$,
\begin{displaymath}
\left\{ \begin{array}{rcll}
N_{0}(t,m) & = & \varphi(\overline{\tau},\pi_{-(t-\overline{\tau})}(m)) K(t-\overline{\tau},m), & \\
N_{n}(t,m) & = & N_{0}(t,m)+G(N_{n-1})(t,m)-J(N_{n-1})(t,m), & n\geq 1.
\end{array} \right.
\end{displaymath}
Alors, $\forall t\geq \overline{\tau}$, $\forall m\in[0,b]$,
\begin{displaymath}
|N_{0}(t,m)|\leq e^{-I(t-\overline{\tau})}||\varphi||_{b}\leq||\varphi||_{b}.
\end{displaymath}
Montrons, par r\'ecurrence, que, pour tout $t\geq \overline{\tau}$ et tout $ m\in[0,b]$,
\begin{displaymath}
 \forall n\in\mathbb{N}, \ |N_{n}(t,m)| \leq ||\varphi||_{b}.
\end{displaymath}
Le r\'esultat est vrai pour $n=0$. Supposons qu'il est vrai au rang $n$, $n\geq0$. Alors, pour tout $t\geq \overline{\tau}$ et tout $ m\in[0,b]$,
\begin{displaymath}
|N_{n+1}(t,m)| \leq |N_{0}(t,m)|+|G(N_{n})(t,m)|+|J(N_{n})(t,m)|.
\end{displaymath}
Or,
\begin{displaymath}
|G(N_{n})(t,m)| \leq 2l(\overline{\tau}-\underline{\tau})\tilde{\zeta}||\varphi||_{b} \int_{\overline{\tau}}^{t} e^{-I(t-s)} ds,
\end{displaymath}
et
\begin{displaymath}
|J(N_{n})(t,m)|\leq l||\varphi||_{b} \int_{\overline{\tau}}^{t} e^{-I(t-s)} ds.
\end{displaymath}
Donc
\begin{displaymath}
|N_{n+1}(t,m)| \leq ||\varphi||_{b}\bigg(e^{-I(t-\overline{\tau})}+l(2(\overline{\tau}-\underline{\tau})\tilde{\zeta}+1)\int_{\overline{\tau}}^{t} e^{-I(t-s)} ds \bigg).
\end{displaymath}
En notant $\mu=\displaystyle\frac{l(2(\overline{\tau}-\underline{\tau})\tilde{\zeta}+1)}{I}$, nous obtenons
\begin{displaymath}
|N_{n+1}(t,m)| \leq ||\varphi||_{b}\bigg(e^{-I(t-\overline{\tau})}+\mu(1- e^{-I(t-\overline{\tau})}) \bigg).
\end{displaymath}
Comme $\mu<1$ par hypoth\`ese, nous aboutissons \`a
\begin{displaymath}
|N_{n+1}(t,m)| \leq ||\varphi||_{b}.
\end{displaymath}
Nous obtenons le r\'esultat attendu.

D'apr\`es la d\'emonstration du th\'eor\`eme \ref{existenceunicite_nonlin}, la suite $(N_{n})_{n\in\mathbb{N}}$ converge vers $N$ lorsque $n$ tends vers $+\infty$. Par passage \`a la limite sur n nous obtenons finalement
\begin{displaymath}
|N^{\varphi}(t,m)|\leq ||\varphi||_{b}, \ \forall t\geq \overline{\tau}, \forall m\in[0,b].
\end{displaymath}
La proposition est ainsi d\'emontr\'ee. $\Box$

\vspace{2ex}
Nous d\'eduisons de ce r\'esultat, le th\'eor\`eme suivant.

\begin{theoreme} \label{invariance}
Si $l(2(\overline{\tau}-\underline{\tau})\tilde{\zeta}+1)<I$ alors, pour tout $b\in(0,g(1)]$ et toute fonction $ \varphi\in\mathcal{C}^{0}([0,\overline{\tau}]\times[0,g(1)])$, il existe $\overline{t} \geq 0$ tel que
\begin{displaymath}
\forall t\geq\overline{t}, \ \ \forall m\in[0,g(1)], \quad |N^{\varphi}(t,m)| \leq ||\varphi||_{b}.
\end{displaymath}
\end{theoreme}
\vspace{1ex}
\textbf{Preuve :} Soient $b\in(0,g(1)]$ et $\varphi\in C([0,\overline{\tau}]\times[0,g(1)])$. Nous d\'efinissons, pour $t\in[0,\overline{\tau}]$, l'application $\tilde{\varphi}$ par :
\begin{displaymath}
\tilde{\varphi}(t,m)= \left\{ \begin{array}{ll}
\varphi(t,m), & \textrm{ si } m\in[0,b], \\
\varphi(t,b), & \textrm{ si } m\in[b,g(1)].
\end{array} \right.
\end{displaymath}
D'apr\`es la proposition \ref{propositioninvariance}, pour tout $t\geq \overline{\tau}$ et tout $m\in[0,g(1)]$, $|N^{\tilde{\varphi}}(t,m)|\leq ||\tilde{\varphi}||_{g(1)}$. De plus, comme $\varphi(t,m)=\tilde{\varphi}(t,m)$ sur $[0,\overline{\tau}]\times[0,b]$, d'apr\`es le th\'eor\`eme \ref{unicite_nonlin} il existe $\overline{t} \geq 0$ tel que $N^{\varphi}(t,m)=N^{\tilde{\varphi}}(t,m)$ pour $t\geq\overline{\tau}$ et $m\in[0,g(1)]$. Enfin, comme $||\tilde{\varphi}||_{g(1)}=||\varphi||_{b}$, nous obtenons
\begin{displaymath}
\forall t\geq\overline{t}, \ \ \forall m\in[0,g(1)], \quad |N^{\varphi}(t,m)|\leq||\varphi||_{b}.
\end{displaymath}
$\Box$

\vspace{2ex}
Ce th\'eor\`eme met en avant le r\'esultat suivant : si l'in\'egalit\'e $l(2(\overline{\tau}-\underline{\tau})\tilde{\zeta}+1)<I$ est v\'erifi\'ee, alors la population de cellules souches domine la population totale, au moins \`a partir d'un certain temps. Or l'in\'egalit\'e $l(2(\overline{\tau}-\underline{\tau})\tilde{\zeta}+1)<I$ est v\'erifi\'ee, par exemple, si $l$ est petit ou $I$ grand, c'est-\`a-dire, biologiquement, si peu de cellules sont introduites de la phase de repos vers la phase de prolif\'eration (car $l\geq\beta(m,x)$ pour tous $m\in[0,g(1)]$ et $x\geq0$) ou si le taux de mortalit\'e dans la phase de repos est \'elev\'e.


\end{document}